\pdfoutput=1
\documentclass[11pt, reqno]{amsart}
\usepackage{amsfonts,latexsym}
\usepackage{amsmath}
\usepackage{amscd}
\usepackage{float,amsmath,amssymb,mathrsfs,bm,multirow,graphics}
\usepackage[dvips]{graphicx}
\usepackage[percent]{overpic}
\usepackage[numbers,sort&compress]{natbib}
\usepackage{enumerate}
\usepackage{amsaddr}

\newcommand{\nequation}{\setcounter{equation}{0}}

\newcommand{\R}{{\Bbb R}}

\newcommand{\C}{{\Bbb C}}

\newcommand{\Z}{{\Bbb Z}}
\newcommand{\proofbegin}{\noindent{\it Proof.\,\,}}
\newcommand{\proofend}{\hfill$\Box$\bigskip}
\newcommand{\proofendcontinue}{\hfill \raisebox{.8mm}[0cm][0cm]{$\bigtriangledown$}\bigskip}

\DeclareMathOperator{\tr}{tr}
\DeclareMathOperator{\sgn}{sgn}
\DeclareMathOperator{\im}{Im}
\DeclareMathOperator{\re}{Re}

\DeclareMathOperator{\sech}{sech}

\def\XXint#1#2#3{{\setbox0=\hbox{$#1{#2#3}{\int}$}
\vcenter{\hbox{$#2#3$}}\kern-.5\wd0}}

\allowdisplaybreaks

\newtheorem{theorem}{Theorem}[section]
\newtheorem{proposition}[theorem]{Proposition}

\newtheorem{definition}[theorem]{Definition}

\newtheorem{remark}[theorem]{Remark}

\newtheorem{figuretext}{Figure}

\input epsf
\title[The nonlinear Schr\"odinger equation with $t$-periodic data I]{\sc The nonlinear Schr\"odinger equation with \\
$t$-periodic data: I. Exact results}

\author{J. Lenells}
\address{Department of Mathematics, KTH Royal Institute of Technology, \\ 100 44 Stockholm, Sweden.}
\email{jlenells@kth.se}

\author{A. S. Fokas}
\address{Department of Applied Mathematics and Theoretical Physics, University of Cambridge, Cambridge CB3 0WA, United Kingdom, and Research Center of Mathematics, Academy of Athens, 11527, Greece.}
\email{T.Fokas@damtp.cam.ac.uk} 
    
\begin{document}

\begin{abstract} 
\noindent
We consider the nonlinear Schr\"odinger equation on the half-line with a given Dirichlet (Neumann) boundary datum which for large $t$ tends to the periodic function $g_0^b(t)$ ($g_1^b(t)$). Assuming that the unknown Neumann (Dirichlet) boundary value tends for large $t$ to a periodic function $g_1^b(t)$ ($g_0^b(t)$), we derive an easily verifiable condition that the functions $g_1^b(t)$ and $g_0^b(t)$ must satisfy. 
Furthermore, we propose two different methods, one based on the formulation of a Riemann-Hilbert problem, and one based on a perturbative approach, for constructing $g_1^b(t)$ ($g_0^b(t)$) in terms of $g_0^b(t)$ ($g_1^b(t)$).
\end{abstract}

\maketitle




\section{Introduction}\nequation
A new method for analyzing boundary value problems for linear and integrable nonlinear PDEs was introduced in \cite{F1997} and \cite{F2000} (see also \cite{F2002, FIS2005, Fbook}) and developed by many authors. For the implementation of this method to integrable nonlinear evolution PDEs see for example \cite{FL2010, trilogy1, trilogy2, Ldnls, LdnlsD2N, LFgnls, LFernst}; reviews for the implementation to linear and to integrable nonlinear PDEs are given in \cite{FS2012, DTV2014} and \cite{F2009, P2014} respectively. 

For integrable nonlinear PDEs this method, which is usually referred to as the {\it unified transform} or the {\it Fokas method}, yields novel integral representations formulated in the complex $k$-plane (the Fourier plane). These representations are similar to the integral representations for the linearized versions of these nonlinear PDEs, but also contain the entries of a certain matrix-valued function, which is the solution of a matrix Riemann-Hilbert (RH) problem. The main advantage of the new method is the fact that this RH problem involves a jump matrix with {\it explicit} $(x,t)$-dependence, uniquely defined in terms of four scalar functions called spectral functions and denoted by $\{a(k), b(k), A(k), B(k)\}$. The functions $a(k)$ and $b(k)$ are defined in terms of the initial datum $u_0(x) = u(x,0)$ via a system of linear Volterra integral equations. The functions $A(k)$ and $B(k)$ are also defined via a system of linear Volterra integral equations, but these integral equations involve {\it all} boundary values. For example, for the nonlinear Schr\"odinger (NLS) equation formulated on the half-line, $A(k)$ and $B(k)$ are defined in terms of the functions $\{u(0,t), u_x(0,t)\}$. However, some of these boundary values are unknown. For example, for the Dirichlet problem of the NLS, the Neumann boundary value $u_x(0,t)$ is unknown. It turns out that this problem can be addressed by employing the so-called {\it global relation}, which is a simple algebraic equation that couples the spectral functions. Actually, by employing this relation, for a particular class of boundary conditions called {\it linearizable}, it is possible to solve the problem on the half-line as effectively as the analogous problem on the full line. Indeed, for linearizable boundary conditions, by utilizing the global relation, it is possible to determine the functions $A(k)$ and $B(k)$ {\it directly}, without the need of determining the unknown boundary values first. 
It should be emphasized that this is true for linearizable boundary conditions of PDEs involving a third order partial $x$-derivative, such as the KdV equation, for which the alternative approach of mapping the half-line problem to a problem on the line apparently fails. 

For non-linearizable boundary conditions, the complete solution of a boundary value problem on the half-line requires the determination of the unknown boundary values, i.e. it requires the characterization of the Dirichlet to Neumann map. This problem was recently analyzed in \cite{trilogy1} and \cite{trilogy2} using two different formulations, both of which are based on the analysis of the global relation: The formulation in \cite{trilogy1} is based on the eigenfunctions involved in the definition of $\{A(k), B(k)\}$ (see also \cite{DMS2001, F2005}), whereas the formulation in \cite{trilogy2} is based on an extension of the Gelfand-Levitan-Marchenko approach first introduced in \cite{BFS2003} and \cite{F2005}.

It must be emphasized that for non-linearizable boundary conditions which decay for large $t$, by utilizing the crucial feature of the new method that it yields RH problems with explicit $(x,t)$-dependence, it is possible to obtain useful asymptotic information about the solution {\it without} characterizing the spectral functions $\{A(k), B(k)\}$ in terms of the given initial and boundary conditions. This can be achieved by employing the Deift-Zhou method \cite{DZ1993} for the long-time asymptotics \cite{FI1992, FI1992b, FI1994, FI1996} and the Deift-Zhou-Venakides method \cite{DVZ1994, DVZ1997} for the zero-dispersion limit \cite{FK2004, K2003}. 

For the physically significant case of boundary conditions which are periodic in $t$, it is not possible to obtain the rigorous form of the long-time asymptotics of the solution, without first characterizing the Dirichlet to Neumann map, at least as $t \to \infty$. Pioneering results in this direction have been obtained for the NLS equation in the quarter plane
\begin{align}\label{nls}
  iu_t + u_{xx} - 2\lambda |u|^2 u = 0, \qquad x>0, \quad t>0, \quad \lambda = \pm 1,
\end{align}
in a series of papers by Boutet de Monvel and coauthors \cite{BIK2007, BIK2009, BKS2009, BKSZ2010} for the particular case that the given Dirichlet datum consists of a {\it single} periodic exponential:
\begin{align}\label{1.1}
  u(0,t) = \alpha e^{i\omega t}, \qquad \alpha > 0, \quad \omega \in \R, \quad t > 0.
\end{align}
In particular, it was shown in \cite{BKS2009} that for the {\it focusing} (i.e. $\lambda = -1$) NLS there exists a solution $u$ which satisfies (\ref{1.1}) as well as the asymptotic condition
\begin{align}\label{1.2}
  u_x(0,t) \sim ce^{i\omega t}, \qquad t\to \infty, \quad c \in \C,
\end{align}
if and only if the triplet of constants $(\alpha,\omega,c)$ satisfy either
\begin{subequations}\label{1.3}
\begin{align}\label{1.3a}
 c = \pm \alpha\sqrt{\omega - \alpha^2} \quad \text{and} \quad \omega \geq \alpha^2,
\end{align}
or
\begin{align}\label{1.3b}
  c = i \alpha \sqrt{|\omega| + 2 \alpha^2} \quad \text{and}  \quad \omega \leq - 6\alpha^2.
\end{align}
\end{subequations}
Equations (\ref{1.3}) show that for $\omega$ in the range $(-6\alpha^2, \alpha^2)$, the asymptotics of $u_x(0,t)$ for a solution satisfying (\ref{1.1}) is not of the simple form $ce^{i\omega t}$. It has been conjectured based on numerical simulations that instead of single exponentials, finite-genus type theta functions arise in this regime, see \cite{BKSZ2010}.

The Dirichlet datum (\ref{1.1}) is complex-valued, thus it cannot be used for the KdV and modified KdV equations in real situations. Results valid for Dirichlet data more general than (\ref{1.1}), including the case of the real-valued Dirichlet datum 
\begin{align}\label{1.4}
  u(0,t) = \alpha \sin{t}, \qquad \alpha \in \R, \quad t > 0,
\end{align}
can be obtained using the perturbative approach introduced in \cite{trilogy2}. In particular, it was shown in \cite{trilogy2} for the NLS, and in \cite{HF2013, H2014} for the mKdV and sine-Gordon equations, that if $u(0,t)$ is given by the right-hand side of (\ref{1.4}), then the function $u_x(0,t)$ for the NLS and the sine-Gordon and the functions $\{u_x(0,t), u_{xx}(0,t)\}$ for the mKdV, respectively, can be computed {\it explicitly} at least up to and including terms of $O(\alpha^3)$, and furthermore the above functions become {\it periodic} as $t \to \infty$. 
Unfortunately, the perturbative approach of \cite{trilogy2} is quite cumbersome and it is practically impossible to go beyond terms of $O(\alpha^3)$.

Here we consider the NLS equation (\ref{nls}) on the half-line and denote by $u_0(x)$ and $g_0(t)$  the given initial datum and the given Dirichlet boundary datum; we also denote by $g_1(t)$ the unknown Neumann boundary value:
\begin{align}\nonumber
& u_0(x) = u(x,0), && 0 < x < \infty; 
	\\ \label{star}
& g_0(t) = u(0,t), \qquad g_1(t) = u_x(0,t), && 0 < t < \infty.
\end{align}
We assume that $u_0 \in \mathcal{S}([0,\infty))$, where $\mathcal{S}([0,\infty))$
denotes the Schwartz class
\begin{align}\label{schwartzdef}
\mathcal{S}([0,\infty)) = \{u \in C^\infty([0,\infty)) \, | \, x^n u^{(m)}(x) \in L^\infty([0,\infty)) \text{ for all } n, m \geq 0\}.
\end{align}
Furthermore, we assume that $g_0(t)$ and $g_1(t)$ are asymptotically periodic as $t \to \infty$, namely,
\begin{align}\label{gjgjB}
g_0(t) - g_0^b(t) = O(t^{-7/2}), \qquad g_1(t) - g_1^b(t) = O(t^{-7/2}), \qquad t\to \infty,
\end{align}
where $g_0^b(t)$ and $g_1^b(t)$ are given periodic function of period $\tau > 0$. 

Boutet de Monvel and coauthors, starting with the particular functions
\begin{align}\label{1.9}
g_0^b(t) = \alpha e^{i\omega t}, \qquad \alpha > 0, \quad \omega \in \R; \qquad
g_1^b(t) = ce^{i\omega t}, \qquad c \in \C,
\end{align}
introduced an ingenious construction involving three steps: In step 1, they introduced an associated background eigenfunction $\psi^b(t,k)$ which they were able to compute explicitly. In step 2, they were able to relate $c$ to $\alpha$ and $\omega$ by the requirement that the relevant spectral functions satisfy the associated global relation. Finally, in step 3, they were able to show that if $g_0^b$  and $g_1^b$ are defined by (\ref{1.9}), then it is possible to define a function $u(x,t)$ via the solution of a $2 \times 2$ matrix RH problem, such that $u(x,t)$  satisfies NLS and $u(0,t)$ and $u_x(0,t)$ asymptote to $g_0^b(t)$ and $g_1^b(t)$, respectively, as $t\to \infty$. 

Here, we show that at least steps 1 and 2 above can be generalized to a large class of $\tau$-periodic functions. The crucial requirement imposed on these functions is that the associated eigenfunction $\psi^b(t,k)$ can be computed explicitly. 
Examples include pairs of functions $\{u(0,t), u_x(0,t)\}$ which can be obtained via the restriction of a function $u(x,t)$ to $x = 0$, where $u(x,t)$ is either a stationary soliton, or a solution obtained via the finite-gap algebro-geometric formalism but which is periodic as opposed to quasiperiodic. 
In general, $\{u(0,t), u_x(0,t)\}$ and hence the associated background eigenfunction $\psi^b(t,k)$ will involve some constants, which are chosen via step 2, i.e. they are chosen by the requirement that the relevant spectral functions satisfy the associated global relation. In this way, we construct pairs of functions $\{g_0^b(t), g_1^b(t)\}$ which we call {\it asymptotically consistent pairs}. For example, for the {\it focusing} NLS (equation (\ref{nls}) with $\lambda = -1$), the pair $\{\alpha e^{i\omega t}, ce^{i\omega t}\}$ is such a pair provided that the triple $(\alpha, \omega, c)$ satisfies (\ref{1.3}).

We emphasize that the requirement for a pair to be asymptotically consistent is simple and easily verifiable. The requirement only involves the zeros of a spectral function $G(k)$ which is related to $\psi^b(t,k)$.

It is an important question to determine which of the asymptotically consistent pairs are actually asymptotically admissible, i.e. for which pairs step 3 can be implemented. 
We propose two different approaches for addressing this question. The first approach can only be applied to the asymptotically consistent pairs which have a linear limit (we say that the pair $\{g_0^b(t), g_1^b(t)\}$ has a linear limit if there exists a solution $u$ of the NLS in the quarter plane whose Dirichlet and Neumann boundary values asymptote to $g_0^b(t)$ and $g_1^b(t)$, respectively, and such that $u = \epsilon u_1 + \epsilon^2 u_2 + \cdots$ satisfies the NLS to each order in $\epsilon$ and each function $u_j(x,t)$ is smooth and has decay as $x \to \infty$).
For example, for $\lambda = -1$ and the single exponential given by the right-hand side of (\ref{1.1}), this approach can be applied to the cases associated with $c = -\alpha\sqrt{\omega - \alpha^2}$ and $c =  i \alpha \sqrt{|\omega| + 2 \alpha^2}$. The approach involves summing up a perturbative series and is presented in detail in \cite{tperiodicII}.
The second approach, which can be applied to any asymptotically consistent pair, involves the investigation of a particular Riemann-Hilbert problem. This approach is a generalization of step 3 mentioned earlier and is discussed in section \ref{consistentsec} of this paper.

\section{Eigenfunctions and spectral functions}\label{eigensec}\nequation 
For simplicity, we will restrict our attention to smooth solutions which decay rapidly as $x \to \infty$.

\begin{definition}\upshape\label{soldef}
A {\it solution of the NLS in the quarter plane} is a smooth function $u:[0,\infty) \times [0,\infty) \to \C$ such that $u(\cdot, t) \in \mathcal{S}([0,\infty))$ for each $t \in [0, \infty)$, and such that (\ref{nls}) is satisfied for $x> 0$ and $t > 0$.
\end{definition}

Let $u(x,t)$ be a solution of the NLS in the quarter plane. We are interested in situations where $g_0(t) = u(0,t)$ and $g_1(t) = u_x(0,t)$ are asymptotically time-periodic. Thus we assume that (\ref{gjgjB}) holds, where $g_j^b(t)$, $j=0,1$, are smooth periodic functions of $t$ with period $\tau>0$.

\subsection{Eigenfunctions}
The NLS equation (\ref{nls}) admits the Lax pair
\begin{align}\label{lax}
\begin{cases}
  \phi_x + ik\sigma_3 \phi = U\phi,
  	\\ 
  \phi_t + 2ik^2 \sigma_3 \phi = V \phi,	
\end{cases}
\end{align}
where $k \in \C$ is the spectral parameter, $\phi(x,t,k)$ is a $2\times 2$-matrix valued eigenfunction, and
\begin{align*}
U = \begin{pmatrix} 0 & u \\
\lambda \bar{u} & 0 \end{pmatrix}, \qquad
V = \begin{pmatrix} -i\lambda |u|^2 & 2ku + iu_x \\
2\lambda k \bar{u} - i\lambda \bar{u}_x & i\lambda |u|^2 \end{pmatrix}, \qquad
\sigma_3 = \begin{pmatrix} 1 & 0 \\ 0 & -1 \end{pmatrix}.
\end{align*}
Let $\hat{\sigma}_3$ act on a $2\times 2$ matrix $A$ by $\hat{\sigma}_3A = [\sigma_3, A]$, i.e. $e^{\hat{\sigma}_3} A = e^{\sigma_3} A e^{-\sigma_3}$.

Following the standard implementation of the unified transform \cite{F1997, F2002}, we define two solutions $\phi_j(x,t,k)$, $j = 2,3$, of (\ref{lax}) by
$$\phi_j(x,t,k) = \mu_j(x,t,k) e^{-i(kx + 2k^2 t)\sigma_3}, \qquad j = 2,3,$$
where $\{\mu_j(x,t,k)\}_2^3$ are the unique solutions of the linear Volterra integral equation
\begin{align}\label{mu23def}
  \mu_j(x,t,k) = &\; I + \int_{(x_j, t_j)}^{(x,t)} e^{i[k(x'-x) + 2k^2(t'-t)] \hat{\sigma}_3} W_j(x',t',k), \qquad j = 2,3,
\end{align}
with $(x_2, t_2) = (0,0)$, $(x_3, t_3) = (\infty, t)$, and
$$W_j = (U dx + V dt)\mu_j, \qquad j =2,3.$$ 
The eigenfunctions $\{\phi_j\}_2^3$ are normalized so that 
$$\phi_2(0,0,k) = I, \qquad \lim_{x \to \infty}(\phi_3(x,t,k) - e^{-i(kx + 2k^2 t)\sigma_3}) = 0.$$

In order to formulate a RH problem suitable for the reconstruction of $u$ from the initial and boundary data, we seek to define an additional eigenfunction $\phi_1(x,t,k)$ which is normalized at $t = \infty$. In the case of decaying boundary data, $\phi_1$ is conveniently defined by $\phi_1 = \mu_1 e^{-i(kx + 2k^2 t)\sigma_3}$, where $\mu_1$ is the solution of (\ref{mu23def}) with $(x_1, t_1) = (0, \infty)$. However, for boundary data satisfying (\ref{gjgjB}), in order to arrive at a well-defined integral equation for $\mu_1$, we first need to subtract off the `background' behavior determined by $\{g_j^b(t)\}_0^1$. In this respect, generalizing the approach of \cite{BK2007, BIK2007, BIK2009, BKS2009}, we define $V^b$ by
$$  V^b(t,k) = \begin{pmatrix} -i\lambda |g_0^b(t)|^2 & 2kg_0^b(t) + ig_1^b(t) \\
2\lambda k \bar{g}_0^b(t) - i\lambda \bar{g}_1^b(t) & i\lambda |g_0^b(t)|^2 \end{pmatrix}, \qquad t > 0, \quad k \in \C.$$
Suppose we can find a solution $\psi^b(t,k)$ of the `background' $t$-part
\begin{align}\label{backgroundtpart}
  \psi_t^b + 2ik^2\sigma_3 \psi^b = V^b \psi^b
\end{align}
of the form 
\begin{align}\label{psiBcalE}
\psi^b(t,k) = \mathcal{E}(t,k) e^{-i\tilde{\Omega}(k) t \sigma_3},
\end{align}
where $\mathcal{E}(t,k)$ is time-periodic with period $\tau$ and $\tilde{\Omega}(k)$ is a complex-valued function.
Then we can define a solution $\phi_1(x,t,k)$ of (\ref{lax}) by
$$\phi_1(x,t,k) = \mu_1(x,t,k) e^{-i(kx + \tilde{\Omega}(k) t)\sigma_3},$$
where $\mu_1(x,t,k)$ is the unique solution of the linear Volterra integral equation
\begin{align}\nonumber
\mu_1(x,t,k) = &\; e^{-ikx\hat{\sigma}_3} \bigg\{\mathcal{E}(t,k) 
-  \mathcal{E}(t,k) \int_t^\infty e^{i\tilde{\Omega}(k) (t' - t)\hat{\sigma}_3} \big[\mathcal{E}^{-1}(t', k) 
	\\\label{mu1def}
&\times (V - V^b)(0,t',k) \mu_1(0,t',k) \big] dt' 
+ \int_0^x e^{ikx' \hat{\sigma}_3}[U(x',t) \mu_1(x',t,k)] dx'\bigg\}.
\end{align}

\subsection{Background eigenfunction}
A background eigenfunction $\psi^b$ of the form (\ref{psiBcalE}) can be constructed using Floquet theory. 
Let $\psi(t,k)$ be the solution of the background $t$-part
\begin{align}\label{tpartB}
  \psi_t + 2ik^2\sigma_3 \psi = V^b\psi,
\end{align}  
normalized by $\psi(0,k) = I$. 

In what follows we will define several quantities in terms of $\psi(t,k)$. We first define the entire $2\times 2$-matrix valued function $Z(k)$ by 
$$Z(k) = \psi(\tau,k).$$ 
The eigenvalues of $Z(k)$ are given by $z(k)$ and $z(k)^{-1}$ where 
\begin{align}\label{lambdaGdef}
z(k) = \frac{1}{2} \big(\tr Z(k) - \sqrt{G(k)}\big), \qquad G(k) = (\tr Z(k))^2 - 4.
\end{align}
Note that $z(k) = 1$ iff $\tr Z(k) =  2$ and $z(k) = -1$ iff $\tr Z(k) = -2$. Also note that $z(k)$ is nonzero for all $k$.

Let $\mathcal{P}$ denote the set of branch points defined by
\begin{align}\label{branchpts}
\mathcal{P} = \{k \in \C \; | \; G(k) = 0,  \; \text{or} \; Z_{12}(k) = 0, \; \text{or} \; Z_{21}(k) = 0\},
\end{align}
where $Z_{ij}$ denotes the $(ij)$'th entry of $Z$. 
The set $\mathcal{P}$ is the union of a finite number of zero sets of entire functions, thus $\mathcal{P}$ is a countable set without accumulation points. Moreover, the symmetries
\begin{align}\label{Zsymm}
Z_{11}(k) = \overline{Z_{22}(\bar{k})}, \qquad Z_{12}(k) = \lambda \overline{Z_{21}(\bar{k})}, 
\end{align}
where $\sigma_1$ is the first Pauli matrix, implies that $\mathcal{P}$ is invariant under the involution $k \mapsto \bar{k}$. 

Let $\mathcal{C}$ denote a set of branch cuts connecting all points in $\mathcal{P}$. We choose these branch cuts so that $\mathcal{C}$ is invariant under the involution $k \mapsto \bar{k}$, see Figure \ref{stationarycuts.pdf} for a possible choice of $\mathcal{C}$ in the case of the stationary one-soliton.
Letting
\begin{align}\label{mathcalSdef}
S^b(k) = \sqrt{-\frac{Z_{11} - Z_{22} - \sqrt{G}}{2\sqrt{G}}} 
\begin{pmatrix} 1 & - \frac{2 Z_{12}}{Z_{11}-Z_{22}-\sqrt{G}} \\
 \frac{2 Z_{21}}{Z_{11}-Z_{22}-\sqrt{G}} & 1 \end{pmatrix}, \qquad k \in \C \setminus \mathcal{C},
 \end{align}
we find that $S^b(k)$ has unit determinant and that
$$Z(k) = S^b(k) \begin{pmatrix} z(k) & 0 \\ 0 & z^{-1}(k) \end{pmatrix} S^b(k)^{-1}, \qquad k \in \C \setminus \mathcal{C}.$$
The identity
$$(Z_{11} - Z_{22} - \sqrt{G})(Z_{11} - Z_{22} + \sqrt{G}) = -4Z_{12}Z_{21},$$
implies that the zeros of $Z_{11} - Z_{22} - \sqrt{G}$ are included in the set of branch points $\mathcal{P}$.

We next define the $2 \times 2$-matrix valued function $\mathcal{B}(k)$ by
\begin{align}\label{Bdef}
\mathcal{B}(k) = \frac{\log z(k)}{\tau} S^b(k) \sigma_3 S^b(k)^{-1}, \qquad k \in \C \setminus \mathcal{C}.
\end{align}
By adding, if necessary, branch cuts to $\mathcal{C}$ to ensure that $\log z(k)$ is single valued on $\C \setminus \mathcal{C}$, 
we find $e^{\tau\mathcal{B}(k)} = Z(k)$ and
$$e^{-t\mathcal{B}(k)} = S^b(k) e^{it\tilde{\Omega}(k) \sigma_3} S^b(k)^{-1}, \qquad k \in \C \setminus \mathcal{C},$$
where 
$$\tilde{\Omega}(k) = -\frac{\log z(k)}{i\tau}.$$ 

The periodicity $V^b(t,k) = V^b(t + \tau,k)$ implies that the matrix valued function $P(t,k)$ defined by
$$P(t,k) = \psi(t,k)  e^{-t\mathcal{B}(k)}, \qquad k \in \C\setminus \mathcal{C},$$
is $t$-periodic with period $\tau$. Indeed, this is a consequence of Floquet theory, see \cite{T2012}.

Defining $\mathcal{E}(t,k)$ by
\begin{align}\label{calEdef}
\mathcal{E}(t,k) = P(t,k) S^b(k), \qquad k \in \C \setminus \mathcal{C},
\end{align}
we deduce that the function $\psi^b(t,k)$ defined by
\begin{align}\label{psiB}
\psi^b(t,k) = \psi(t,k) S^b(k)
= \mathcal{E}(t,k) e^{-i\tilde{\Omega}(k) t \sigma_3}, \qquad k \in \C \setminus \mathcal{C},
\end{align}
is a solution of (\ref{tpartB}) of the form (\ref{psiBcalE}). Note that $\psi^b(t,k)$ is periodic iff $z(k) = 1$ and antiperiodic iff $z(k) = -1$. In this sense, the zero set of $G(k)$ is the union of the periodic and antiperiodic spectrum.

\subsection{Asymptotics as $k \to \infty$}\label{kinftysubsec}
As $k \to \infty$, we have
$$\psi(t,k) = \bigg(I + \frac{\psi^{(1)}(t)}{k} + \frac{\psi^{(2)}(t)}{k^2} + \cdots \bigg)e^{-2ik^2 t \sigma_3}
+ \bigg(\frac{\tilde{\psi}^{(1)}(t)}{k} + \frac{\tilde{\psi}^{(2)}(t)}{k^2} + \cdots \bigg)e^{2ik^2 t \sigma_3},$$
where the coefficient matrices $\psi^{(j)}(t)$ and $\tilde{\psi}^{(j)}(t)$ are independent of $k$ and $\tilde{\psi}^{(1)}(t)$ is off-diagonal. The first few coefficients in this expansion are derived by integration by parts and are given explicitly by
\begin{align}\nonumber
[\psi(t,k)]_2 = &\; \Bigg\{\begin{pmatrix} 0 \\ 1 \end{pmatrix} + \frac{1}{k}\begin{pmatrix} - \frac{ig_0^b(t)}{2} \\ i\eta_1(t) \end{pmatrix}
 + \frac{1}{k^2} \begin{pmatrix} \frac{g_1^b(t)}{4} + \frac{g_0^b(t)}{2}\eta_1(t) \\
i\eta_2(t) + \frac{\lambda}{4}|g_0^b(0)|^2
\end{pmatrix} + O\bigg(\frac{1}{k^3}\bigg)\Bigg\}e^{2ik^2t}
	\\ \label{psi2expansion}
& + \Bigg\{ \frac{1}{k}\begin{pmatrix} \frac{ig_0^b(0)}{2} \\ 0 \end{pmatrix}
 + \frac{1}{k^2} \begin{pmatrix} \frac{g_0^b(0)}{2}\eta_1(t) - \frac{g_1^b(0)}{4} \\ - \frac{\lambda}{4} \overline{g_0^b(t)} g_0^b(0) \end{pmatrix}  + O\bigg(\frac{1}{k^3}\bigg)\Bigg\} e^{-2ik^2t},
\end{align}
where $[M]_1$ and $[M]_2$ denote the first and second columns of a $2 \times 2$ matrix $M$ and the real-valued functions $\{\eta_j(t)\}_1^2$ are defined by
\begin{align*}
& \eta_1(t) = \lambda \int_0^t \im(\overline{g_0^b(t')}g_1^b(t')) dt', 
	\\
& \eta_2(t) = 
\frac{\lambda}{4} \int_0^t \big[\lambda |g_0^b|^4 - |g_1^b|^2 - 4i\im(\bar{g}_0^bg_1^b)\eta_1 - i\bar{g}_0^bg_{0t}^b\big] dt'.
\end{align*}
Evaluating (\ref{psi2expansion}) at $t = \tau$ and using the periodicity of $g_0^b$ and $g_1^b$, we find
\begin{align*}
[Z(k)]_2 = &\; \Bigg\{\begin{pmatrix} 0 \\ 1 \end{pmatrix} + \frac{1}{k}\begin{pmatrix} - \frac{ig_0^b(0)}{2} \\ i\eta_1(\tau) \end{pmatrix}
 + \frac{1}{k^2} \begin{pmatrix} \frac{g_1^b(0)}{4} + \frac{g_0^b(0)}{2}\eta_1(\tau) \\
i\eta_2(\tau) + \frac{\lambda}{4}|g_0^b(0)|^2
\end{pmatrix}  + O\bigg(\frac{1}{k^3}\bigg)\Bigg\}e^{2ik^2\tau}
	\\
& + \Bigg\{ \frac{1}{k}\begin{pmatrix} \frac{ig_0^b(0)}{2} \\ 0 \end{pmatrix}
 + \frac{1}{k^2} \begin{pmatrix} \frac{g_0^b(0)}{2} \eta_1(\tau)- \frac{g_1^b(0)}{4} \\ - \frac{\lambda}{4} |g_0^b(0)|^2 \end{pmatrix}  + O\bigg(\frac{1}{k^3}\bigg)\Bigg\} e^{-2ik^2\tau}.
\end{align*}
Hence,
\begin{align}\nonumber
& Z(k) = \begin{pmatrix} e^{-2ik^2\tau} & 0 \\ 0 & e^{2ik^2\tau} \end{pmatrix} 
+ \frac{1}{k}\begin{pmatrix} -i\eta_1(\tau)e^{-2ik^2\tau} & g_0^b(0)\sin(2k^2\tau) \\ 
\lambda \overline{g_0^b(0)}\sin(2k^2\tau) & i\eta_1(\tau)e^{2ik^2\tau} \end{pmatrix}
	\\ \nonumber
& + \frac{1}{k^2} \begin{pmatrix} -i\overline{\eta_2(\tau)}e^{-2ik^2\tau}  -\frac{i\lambda}{2}|g_0^b(0)|^2 \sin(2k^2\tau) & \frac{ig_1^b(0)\sin(2k^2\tau)}{2} + g_0^b(0)\eta_1(\tau)\cos(2k^2\tau) \\ 
-\frac{i\overline{g_1^b(0)}\sin(2k^2\tau)}{2} + g_0^b(0)\eta_1(\tau)\cos(2k^2\tau) 
& i\eta_2(\tau)e^{2ik^2\tau} + \frac{i\lambda}{2}|g_0^b(0)|^2 \sin(2k^2\tau)
\end{pmatrix}
	\\\label{Zasymptotics}
& + O\bigg(\frac{e^{2ik^2\tau}}{k^3}\bigg) + O\bigg(\frac{e^{-2ik^2\tau}}{k^3}\bigg), \qquad k \to \infty, \quad k \in \C.
\end{align}
In particular,
\begin{align}\nonumber
G(k) = & -4 \sin^2(2k^2 \tau)
-\frac{8\eta_1(\tau)}{k} \cos(2k^2 \tau)\sin(2k^2 \tau)+ O\bigg(\frac{1}{k^2}\bigg) + O\bigg(\frac{e^{4ik^2\tau}}{k^2}\bigg) 
	\\ \label{Gasymptotics}
& + O\bigg(\frac{e^{-4ik^2\tau}}{k^2}\bigg),\qquad k \to \infty, \quad k \in \C.
\end{align}

\subsection{Choice of branches}
We fix the branches of $\sqrt{G}$, $\sqrt{-\frac{Z_{11} - Z_{22} - \sqrt{G}}{2\sqrt{G}}}$, and $\log z(k)$, by requiring that 
\begin{align}\nonumber
& \sqrt{G(k)} = 2i\sin(2k^2 \tau)(1 + O(k^{-1})), &&
\sqrt{-\frac{Z_{11} - Z_{22} - \sqrt{G}}{2\sqrt{G}}} = 1 + O(k^{-1}),
	\\ \label{branchchoices}
& \log z(k) = -2ik^2 \tau + O(k^{-1}), && \tilde{\Omega}(k) = 2k^2 + O(k^{-1}),
\end{align}
as $k$ goes to infinity in $\C$ with $k$ remaining a bounded distance away from $\mathcal{C} \cup 
\{\text{zeros of $\sin(2k^2 \tau)$}\} \subset \mathcal{C} \cup \R \cup i\R$.

\subsection{Boundedness and analyticity properties}
Let 
\begin{align}\nonumber
D_1 = \{\im k > 0\} \cap \{\im \tilde{\Omega}(k) > 0\},  \qquad
D_2 = \{\im k > 0\} \cap \{\im \tilde{\Omega}(k) < 0\}, 
	\\ \label{Djdef}
D_3 = \{\im k < 0\} \cap \{\im \tilde{\Omega}(k) > 0\},  \qquad
D_4 = \{\im k < 0\} \cap \{\im \tilde{\Omega}(k) < 0\},
\end{align}
and let $D_+ = D_1 \cup D_3$ and $D_- = D_2 \cup D_4$. The asymptotics (\ref{branchchoices}) of $\tilde{\Omega}(k)$ implies that $D_j$ can be viewed as a deformation of the $j$'th quadrant. We will assume that the branch cuts $\mathcal{C}$ are chosen in such a way that $D_1 \setminus \mathcal{C}$ is connected. 

Henceforth $C > 0$ will denote a generic constant.

\begin{proposition}
The eigenfunctions $\{\mu_j(x,t,k)\}_1^3$ defined by (\ref{mu23def}) and (\ref{mu1def}) possess the following analyticity and boundedness properties:
\begin{enumerate}[$(a)$]
  \item The first (resp. second) column of $\mu_1(0,t,k)$ is defined and analytic for $D_- \setminus \mathcal{C}$ (resp. $D_+ \setminus \mathcal{C}$).
The second column of $\mu_1$ has a continuous extension to the boundary of $D_+ \setminus \mathcal{C}$ in the sense that, away from the branch points $\mathcal{P}$, the limits as $k$ approaches the boundary of $D_+ \setminus \mathcal{C}$ exist and the boundary function is continuous. If the boundary contains a branch cut which can be approached from both left and right within $D_+ \setminus \mathcal{C}$, then the left and right limits are, in general, different.
   
  \item $\mu_1(0,t,k)$ approaches $\mathcal{E}(t,k)$ as $t \to \infty$. More precisely, if $K_\pm$ are compact subsets of $(\overline{D_\pm \setminus \mathcal{C}})\setminus \mathcal{P}$, then
\begin{align}\label{mu1calE}
 |\mu_1(0,t,k) - \mathcal{E}(t,k)| \leq C(1+t)^{-5/2}, \qquad k \in (K_-, K_+), \quad t \in [0, \infty).
\end{align}
 
  \item $\mu_2(x,t,k)$ is defined and analytic for all $k \in \C$.

  \item The first (resp. second) column of $\mu_3(x,t,k)$ is defined and analytic for $\im k < 0$ (resp. $\im k > 0$) with a continuous extension to $\im k \leq 0$ (resp. $\im k \geq 0$).

    \item The $\mu_j$'s have unit determinant whenever the determinant is defined. In particular,
\begin{align*}
& \det \mu_2(x,t,k) = 1, \qquad k \in \C,
	\\
& \det \mu_3(x,t,k) = 1, \qquad k \in \R.
\end{align*}
\end{enumerate}
\end{proposition}
\proofbegin
We will prove $(a)$ and $(b)$; the proofs of $(c)$-$(e)$ are standard \cite{F2002, FIS2005}. 

The second column of (\ref{mu1def}) evaluated at $x = 0$ can be written as
\begin{align}\label{integraleq}
\Psi(t,k) = [\mathcal{E}(t,k)]_2 - \int_t^\infty  E_1(t,t',k) \Delta(t',k) \Psi(t',k) dt',
\end{align}
where $\Psi(t,k) = [\mu_1(0,t,k)]_2$, $\Delta(t',k) = (V - V^b)(0,t',k)$, and
$$E_1(t,t',k) = \mathcal{E}(t,k)\begin{pmatrix} e^{2i\tilde{\Omega}(k)(t'-t)} & 0 \\ 0 & 1 \end{pmatrix} \mathcal{E}^{-1}(t',k).$$
Let $F$ denote the closure of $D_+ \setminus \mathcal{C}$ with the branch points $\mathcal{P}$ removed, i.e. $F = (\overline{D_+ \setminus \mathcal{C}})\setminus \mathcal{P}$. Suppose $K$ is a compact subset of $F$. 
Let $\Psi_0 = [\mathcal{E}(t,k)]_2$ and define $\Psi_l$ for $l \geq 1$ inductively by 
\begin{align*}
\Psi_{l+1}(t,k) = - \int_t^\infty E_1(t, t',k) \Delta(t', k) \Psi_l(t', k) dt', \qquad t \in \R, \quad k \in K.
\end{align*}  
Then
\begin{align}\label{Philiterated}
\Psi_l(t,k) = & (-1)^l \int_{t = t_{l+1} \leq t_l \leq \cdots \leq t_1 < \infty} \prod_{i = 1}^l E_1(t_{i+1}, t_i, k) \Delta(t_i, k) [\mathcal{E}(t_1,k)]_2 dt_1 \cdots dt_l. 
\end{align}
Using the estimate 
\begin{align*}
  |E_1(t, t', k)| < C, \qquad 0 \leq t \leq t' < \infty, \quad k \in K,
\end{align*}
we find, for $k \in K$ and $t \in [0, \infty)$,
\begin{align}\nonumber
|\Psi_l(t,k)| \leq \; & C \int_{t \leq t_l \leq \cdots \leq t_1 < \infty} \prod_{i = 1}^l  |\Delta(t_i, k)| |[\mathcal{E}(t_1,k)]_2| dt_1 \cdots dt_l 
	\\ \label{Philestimate}
\leq & \; \frac{C}{l!}\|[\mathcal{E}(\cdot,k)]_2\|_{L^\infty([t,\infty))} \|\Delta(\cdot, k)\|_{L^1([t,\infty))}^l.
\end{align}
Since, by periodicity of $\mathcal{E}(t,k)$ and (\ref{gjgjB}), 
$$\sup_{k \in K} \|[\mathcal{E}(\cdot,k)]_2\|_{L^\infty([0,\infty))} < C, \qquad
\|\Delta(\cdot, k)\|_{L^1([t,\infty))} < C(1+t)^{-5/2},$$ 
we find
$$|\Psi_l(t,k)| \leq \frac{C}{l!}((1+t)^{-5/2})^l, \qquad t \in [0, \infty), \quad k \in K.$$
Hence the series
$$\Psi(t,k) = \sum_{l=0}^\infty \Psi_l(t,k)$$
converges absolutely and uniformly for $t \in [0, \infty)$ and $k \in K$ to a continuous solution $\Psi(t,k)$ of (\ref{integraleq}) which satisfies
\begin{align}\label{Psiminusf}
|\Psi(t,k) - [\mathcal{E}(t,k)]_2| \leq \sum_{l=1}^\infty | \Psi_l(t,k)| \leq C (1+t)^{-5/2}, \qquad t \in [0, \infty), \quad k \in K.
\end{align}
Since $\tilde{\Omega}$ is analytic in $\C \setminus \mathcal{C}$, equation (\ref{Philiterated}) together with (\ref{gjgjB}) and an easy estimate of $|\partial_k E_1(t, t', k)|$ shows that $\Psi_l$ is analytic in $K \cap (D_+ \setminus \mathcal{C})$. The analyticity of $\Psi(t,k)$ follows from the uniform convergence.
This establishes $(a)$ and $(b)$ for the second column of $\mu_1$; the proof for the first column is similar. 
\proofend

\subsection{Spectral functions}
We define the spectral functions $\{a(k),b(k), A(k),B(k)\}$ by
\begin{align}\label{sSdef}
s(k) = \begin{pmatrix}\overline{a(\bar{k})} & b(k) \\
  \lambda \overline{b(\bar{k})} & a(k) \end{pmatrix}, \qquad
  S(k) = \begin{pmatrix}\overline{A(\bar{k})} & B(k) \\
  \lambda \overline{B(\bar{k})} & A(k) \end{pmatrix},
  \end{align}
where  
\begin{align}\label{sS}
s(k) = \mu_3(0,0,k), 
\qquad S(k) = \mu_1(0,0,k) = \lim_{t\to \infty} e^{2ik^2t\sigma_3}\mu_2^{-1}(0,t,k) \mathcal{E}(t,k) e^{-i\tilde{\Omega}(k)t \sigma_3}.
\end{align}
The function $\mu_3(x,0,k)$, $x \geq 0$, and hence also $s(k)$, is defined in terms of  the initial datum $u_0(x)$. The function $\mu_1(0,t,k)$, $t \geq 0$, and hence also $S(k)$, is defined in terms of $\{g_0(t), g_1(t), g_0^b(t), g_1^b(t)\}$. 
The functions $a(k)$ and $b(k)$ are defined and analytic in $\im k > 0$ with a continuous extension to $\im k \geq 0$. The functions $A(k)$ and $B(k)$ are defined and analytic in $D_+ \setminus \mathcal{C}$ and have continuous limits as $k$ approaches the boundary of this set away from the branch points.

In analogy with (\ref{sSdef}) we define $\{A^b(k), B^b(k)\}$ by
$$S^b(k) = \begin{pmatrix} \overline{A^b(\bar{k})} & B^b(k) \\ \lambda\overline{B^b(\bar{k})} & A^b(k) \end{pmatrix}.$$

\subsection{Symmetries}
The symmetry $Z(k) = \sigma_1\overline{Z(\bar{k})} \sigma_1$ implies that
\begin{align*}
 G(k) = \overline{G(\bar{k})}, \qquad \sqrt{G(k)} = -\overline{\sqrt{G(\bar{k})}}, \qquad
 z(k) = \frac{1}{\overline{z(\bar{k})}},\qquad
\tilde{\Omega}(k) = \overline{\tilde{\Omega}(\bar{k})},
\end{align*}
for $k \in \C \setminus \mathcal{C}$. Moreover, by (\ref{mathcalSdef}),
$$A^b(k) = \overline{A^b(\bar{k})}, \qquad k \in \C \setminus \mathcal{C}.$$

\subsection{The global relation}
We henceforth make the following rather mild assumption:
\begin{align*}
\begin{cases}
 \text{There exists a $c > 0$ such that $[\mu_3(0,t,k)]_2 = O(e^{ct})$}
	\\
 \text{as $t \to \infty$, uniformly for all sufficiently large $k \in D_1$.}
 \end{cases}
\end{align*}
This assumption holds, for example, if $\|u(\cdot, t)\|_{L^1([0,\infty))} = O(t)$ as $t \to \infty$.
As an example, we note that for the stationary soliton of section \ref{stationarysec},  $\|u(\cdot, t)\|_{L^1([0,\infty))} = 2 \arctan(e^\gamma)$ is a constant independent of $t$ and $[\mu_3(0,t,k)]_2 = O(1)$ as $t \to \infty$ uniformly for all $k \in \C$ bounded away from $-i\sqrt{\omega}/2$.

Letting $t \to \infty$ in the $(12)$ entry of the relation
\begin{align}\label{Sinvs}
S^{-1}(k)s(k) = e^{i\tilde{\Omega}(k)t\sigma_3} \mu_1^{-1}(0,t,k) \mu_3(0,t,k) e^{-2ik^2t\sigma_3}
\end{align}
and using the above assumption, we find the following global relation:
\begin{align}\label{GR}
A(k) b(k) - B(k) a(k) = 0, \qquad  k \in D_1 \setminus \mathcal{C}, \quad \im(\tilde{\Omega}(k) + 2k^2) > c, 
\end{align}
where the condition $k \in D_1 \setminus \mathcal{C}$ ensures that the left-hand side is well-defined, while the condition $\im(\tilde{\Omega}(k) + 2k^2) > c$ ensures that the right-hand side of (\ref{Sinvs}) vanishes as $t \to \infty$. In view of (\ref{branchchoices}), the condition $\im(\tilde{\Omega}(k) + 2k^2) > c$ is fulfilled for all sufficiently large $k \in D_1 \setminus \mathcal{C}$ in a sector $\arg k \in (\epsilon, \frac{\pi}{2} - \epsilon)$. Since we are assuming that the branch cuts are chosen in such a way that $D_1 \setminus \mathcal{C}$ is connected, analytic continuation implies that (\ref{GR}) is valid for all $k \in \overline{D_1 \setminus \mathcal{C}}$.

With a slight abuse of notation, we write
\begin{align*}
\mu_3(0,t,k) = \begin{pmatrix}\overline{a(t,\bar{k})} & b(t,k) \\
  \lambda \overline{b(t,\bar{k})} & a(t,k) \end{pmatrix}, \qquad
  \mu_1(0,t,k) = \begin{pmatrix}\overline{A(t,\bar{k})} & B(t,k) \\
  \lambda \overline{B(t,\bar{k})} & A(t,k) \end{pmatrix}.
\end{align*}
The expression
$$C(t, k) = e^{i\tilde{\Omega}(k)t\sigma_3}\mu_1^{-1}(0,t,k) \mu_3(0,t,k) e^{-2ik^2t\sigma_3}$$
is independent of $t$ as a consequence of (\ref{Sinvs}).
Letting $t' \to \infty$ in the $(12)$ entry of the relation
\begin{align}
C(t, k) = C(t',k),
\end{align}
we find the following generalization of (\ref{GR}):
\begin{align}\label{gGR}
& A(t,k) b(t,k) - B(t,k) a(t,k) = 0, \qquad t \geq 0, \quad k \in \overline{D_1 \setminus \mathcal{C}}.
\end{align}

\section{Asymptotically admissible pairs}\nequation

We next introduce the notion of an asymptotically admissible pair.

\begin{definition}\upshape\label{asymptoticallyadmissibledef}
A pair of smooth functions $\{g_0^b(t), g_1^b(t)\}$, $t \geq 0$, is {\it asymptotically admissible} for NLS if there exists a solution $u(x,t)$ of the NLS in the quarter plane (see Definition \ref{soldef}) such that the Dirichlet and Neumann boundary values of $u$ asymptote towards $g_0^b(t)$ and $g_1^b(t)$ respectively in the sense that
\begin{align}\label{u0uxasymptote}
u(0, t) - g_0^b(t) = O(t^{-7/2}), \qquad u_{x}(0,t) - g_1^b(t) = O(t^{-7/2}), \qquad t\to \infty.
\end{align}
A pair which is not asymptotically admissible is called {\it asymptotically inadmissible}.
\end{definition}

The following theorem gives a necessary condition for a periodic pair to be asymptotically admissible.

\begin{theorem}\label{asymptoticallyadmissibleth}
 Let $\{g_0^b(t), g_1^b(t)\}$ be smooth periodic functions of period $\tau > 0$. If the pair $\{g_0^b(t), g_1^b(t)\}$ is asymptotically admissible for NLS, then the function $G(k)$, defined in (\ref{lambdaGdef}), has no zeros of odd order in the interior of $\bar{D}_1 \cup \bar{D}_4$.
\end{theorem}
\proofbegin
Suppose $\{g_0^b(t), g_1^b(t)\}$ is asymptotically admissible and let $u(x,t)$ be a solution of the defocusing NLS in the quarter plane satisfying (\ref{u0uxasymptote}). 

\medskip\noindent
{\bf Claim 1.} $G(k)$ has no zeros of odd order in the interior of $\bar{D}_1$.

{\it Proof of Claim 1.}
Suppose $\kappa$ is such a zero of $G(k)$. Then $z(\kappa) = 1$ or $z(\kappa) = -1$, so $\im \tilde{\Omega}(\kappa) = \frac{1}{\tau}\log |z(\kappa)| = 0$. Now $z = \frac{1}{2}(\tr Z - \sqrt{G})$  is an analytic map $U \to \C$, where $U$ denotes a sufficiently small neighborhood of $\kappa$ in the Riemann surface defined by 
\begin{align}\label{RiemannSurface}
\{(k,l) \in \C^2 \, | \, l^2 = k - \kappa\}.
 \end{align}
Indeed, if $n \in \Z$ denotes the odd order of the zero $\kappa$, then
\begin{align}\label{sqrtGG2l2}
\sqrt{G(k)} = \sqrt{(k - \kappa)^n G_1(k - \kappa)} = (k - \kappa)^{n/2} \sqrt{G_1(l^2)} = l^n G_2(l^2)
\end{align}
where $G_j$ is analytic near $\kappa$ and $G_j(0) \neq 0$ for $j = 1,2$.
Equation (\ref{sqrtGG2l2}) also shows that if we let $k^+ = (k,l)$ and $k^- = (k,-l)$ denote the points in the upper and lower sheets of $U$ lying over $k$, then $\sqrt{G(k^+)} = - \sqrt{G(k^-)}$. Hence $z(k^+) = z(k^-)^{-1}$ and $\im \tilde{\Omega}(k^+) = -\im \tilde{\Omega}(k^-)$.
Shrinking $U$ if necessary, $\log z(\cdot)$ is an analytic function  from $U$ to $\R$. It follows that there exists a curve $\gamma$ in $U$  passing through $\kappa$ such that $\im \tilde{\Omega} =  \frac{1}{\tau} \re \log z = 0$ on $\gamma$. Since $\im \tilde{\Omega}(k^+) = 0$ iff $\im \tilde{\Omega}(k^-) = 0$, we may assume that $\gamma$  is invariant under the sheet changing involution $(k,l) \to (k,-l)$. 

We deform the branch cut $C$ for $\sqrt{G}$ that begins at $\kappa$  so that it coincides with the natural projection of $\gamma$ onto $\C$ for $k \in B_r$, where $B_r \subset \C$ denotes an open disk of radius $r > 0$ centered at $\kappa$ which contains no other branch points and which is contained in $\bar{D}_1$. Then $\im \tilde{\Omega} = 0$ on $C \cap B_r$.
Let $\mathcal{E}_+$ and $\mathcal{E}_-$ denote the limits of $\mathcal{E}$ onto $C \cap B_r$ from the left and right, respectively. Let $(\mu_1(0,t,k))_{\pm}$ be defined in terms of $\mathcal{E}_\pm(t,k)$ via (\ref{mu1def}).
Since $\im \tilde{\Omega}(k) = 0$ for $k \in C \cap B_r$, both columns of $(\mu_1)_\pm$ are well-defined. Define the functions $\nu_\pm(t,k)$ by
\begin{align}\label{mu1nuE}
  (\mu_1(0,t,k))_\pm = \nu_\pm(t,k) \mathcal{E}_\pm(t,k), \qquad k \in C \cap B_r.
\end{align}
Then $\nu_\pm$ satisfy the integral equation
\begin{align}\nonumber
\nu_\pm(t,k) = &\; I - \int_t^\infty \psi^b(t,k) (\psi^b)^{-1}(t',k) (V - V^b)(0,t', k) \nu_\pm(t', k)
	\\\label{nuVolterra}
&\times \psi^b(t', k) (\psi^b)^{-1}(t, k) dt'.
\end{align}
By (\ref{psiB}),
\begin{align}\label{psiBpsiBpsipsi}
\psi^b(t,k)(\psi^b)^{-1}(t',k) = \psi(t,k)\psi^{-1}(t',k).
\end{align}
Hence $\psi^b(t,k)(\psi^b)^{-1}(t',k)$ and its inverse are entire functions of $k$. 
Moreover, the identity
\begin{align}\label{psibpsibidentity}
\psi^b(t,k)(\psi^b)^{-1}(t',k) = \mathcal{E}(t,k) e^{-i\tilde{\Omega}(k)(t-t')\sigma_3} \mathcal{E}^{-1}(t',k),
\end{align}
which also follows from (\ref{psiB}), shows that $\psi^b(t,k)(\psi^b)^{-1}(t',k)$ and its inverse are bounded for $t, t' \geq 0$ whenever $\im \tilde{\Omega}(k) = 0$. The assumption $V - V^b = O(t^{-7/2})$ then implies that the Volterra equation (\ref{nuVolterra}) has a unique solution for $k \in C \cap B_r$. Hence $\nu_- = \nu_+$.

Evaluating the second column of (\ref{mu1nuE}) at $t = 0$ and using the fact that $\mathcal{E}(0,k) = S^b(k)$, we find
\begin{align}\label{BApm}
\begin{pmatrix} B(k) \\ A(k) \end{pmatrix}_\pm = \nu(0,k) \begin{pmatrix} B^b(k) \\ A^b(k)\end{pmatrix}_\pm.
\end{align}
Using the short-hand notation $\nu_{ij}$ for the $(ij)$'th entry of $\nu(0,k)$, equations (\ref{mathcalSdef}) and (\ref{BApm}) imply
\begin{align}\label{fracBA}
\bigg(\frac{B(k)}{A(k)}\bigg)_\pm = \bigg(\frac{\nu_{11}B^b + \nu_{12}A^b}{\nu_{21}B^b + \nu_{22}A^b}\bigg)_\pm, \qquad k \in C \cap B_r.
\end{align}
It follows that
\begin{align}\nonumber
& \bigg(\frac{B(k)}{A(k)}\bigg)_+ - \bigg(\frac{B(k)}{A(k)}\bigg)_-
	\\ \label{BoverAjump}
& = \frac{2 Z_{12} ((\sqrt{G})_--(\sqrt{G})_+)}{(\nu_{22}  ((\sqrt{G})_- - Z_{11}+Z_{22})+2 \nu_{21} Z_{12}) (\nu_{22} ((\sqrt{G})_+-Z_{11}+Z_{22})+2 \nu_{21} Z_{12})}.
\end{align}
Since $G$ does not vanish identically and $\kappa$ is a zero of odd order, the jump $(\sqrt{G})_--(\sqrt{G})_+ = 2(\sqrt{G})_-$ is nonzero across $C$. Equation (\ref{BoverAjump}) therefore implies that the quotient $B(k)/A(k)$ is discontinuous in $\bar{D}_1$. Since $a(k)$ and $b(k)$ are continuous in $\bar{D}_1 \cup \bar{D}_2$, this contradicts the global relation (\ref{GR}).
\proofendcontinue

\medskip\noindent
{\bf Claim 2.}  $G(k)$ has no zeros of odd order in the interior of $\bar{D}_4$.

{\it Proof of Claim 2.}
Since $G(k) = \overline{G(\bar{k})}$, this is a direct consequence of Claim 1 and the fact that the involution $k \mapsto \bar{k}$ carries $D_1$ onto $D_4$.
\proofendcontinue

It only remains to prove that $G(k)$ has no zeros of odd order in $\text{Int}(\bar{D}_1 \cup \bar{D}_4) \cap \R$. 
We treat the focusing and defocusing cases separately. Note that $G(k) \in \R$ for $k \in \R$.

\medskip\noindent
{\bf Claim 3.} For the focusing NLS (i.e. $\lambda = -1$), we have $G(k) \leq 0$ for all $k \in \R$. In particular, $G(k)$ has no zeros of odd order in $\R$.

{\it Proof of Claim 3.}
The function $Z(k)$ is entire and 
$$Z_{11}(k) = \overline{Z_{22}(k)}, \quad Z_{12}(k) = -\overline{Z_{21}(k)}, \qquad k \in \R.$$
Therefore the determinant relation $\det Z(k) = 1$ yields $|Z_{11}|^2 = 1 - |Z_{12}|^2 \leq 1$  on $\R$. Hence
$$G(k) = 4((\re Z_{11}(k))^2 - 1) \leq 4(|Z_{11}(k)|^2 - 1) \leq 0, \qquad k \in \R.$$
\proofendcontinue

For the remainder of the proof we consider the defocusing NLS, i.e., we assume $\lambda = 1$. 

Define $Q^b(k)$ by
\begin{align}\label{fdef}  
  Q^b(k) = \frac{B^b(k)}{A^b(k)} = - \frac{2 Z_{12}(k)}{Z_{11}(k) -Z_{22}(k) -\sqrt{G(k)}}, \qquad k \in \C \setminus \mathcal{C}.
  \end{align}
For $\kappa \in \R$, we let $Q^b_+(\kappa)$ and $Q^b_-(\kappa)$ denote the limits of $Q^b(k)$ as $k$ approaches $\kappa$  from the upper and lower half-planes, respectively.
 
\medskip\noindent
{\bf Claim 4.} If $G(k) > 0$ for some $k \in \R$, then $|Q^b_\pm(k)| = 1$.

{\it Proof of Claim 4.}
 If $G(k) = 4((\re Z_{11}(k))^2 - 1) > 0$, then $|\re Z_{11}(k)| > 1$; thus
$$Q^b_\pm(k) 
= \frac{-Z_{12}(k)}{i\im Z_{11}(k) + \sigma \sqrt{(\re Z_{11}(k))^2 - 1}},$$
where $\sigma = 1$  or $\sigma = -1$.
The determinant condition $\det Z(k) = 1$ yields $|Z_{11}|^2 - |Z_{12}|^2 = 1$  on $\R$. Hence
$$|Q^b_\pm(k)| = \frac{|Z_{12}(k)|}{\sqrt{(\im Z_{11}(k))^2 + (\re Z_{11}(k))^2 - 1}} = 1.$$
\proofendcontinue

Now suppose $\kappa \in \text{Int}(\bar{D}_1 \cup \bar{D}_4) \cap \R$ is a zero of odd order of $G(k)$.
Then $G(k)$ changes sign at $\kappa$. For definiteness, we suppose $G(k) > 0$ to the left of $\kappa$. Let $[\kappa - \epsilon, \kappa + \epsilon]$ denote a closed interval centered at $\kappa$ which contains no other zeros or branch cuts of $\sqrt{G(k)}$ and such that $[\kappa - \epsilon, \kappa + \epsilon] \subset \bar{D}_1$.
We choose the branch cut that begins at $\kappa$ so that its restriction to a ball of radius $\epsilon$ centered at $\kappa$ coincides with $[\kappa - \epsilon, \kappa]$.

The relation $z(k) \overline{z(\bar{k})} = 1$ and the continuity of $z(k)$ across $[\kappa, \kappa + \epsilon]$ imply that $|z(k)|^2 = 1$ for $k \in [\kappa, \kappa + \epsilon]$. Hence $\im \tilde{\Omega} = 0$ on $[\kappa, \kappa + \epsilon]$, so we may define $\nu(t,k)$ for $k \in [\kappa, \kappa + \epsilon]$ by 
\begin{align}\nonumber
\nu(t,k) = &\; I - \int_t^\infty \psi^b(t,k) (\psi^b)^{-1}(t',k) (V - V^b)(0,t', k) \nu(t', k)
	\\\label{nuVolterra2}
&\times \psi^b(t', k) (\psi^b)^{-1}(t, k) dt', \qquad t \geq 0, \quad k \in [\kappa, \kappa + \epsilon].
\end{align}

\medskip\noindent
{\bf Claim 5.} The function $\nu(t,k)$ defined by (\ref{nuVolterra2}) is a continuous function of $(t,k) \in [0, \infty) \times [\kappa, \kappa + \epsilon]$. 
In particular, $\nu(0,k)$ is well-defined, continuous, and bounded for $k \in [\kappa, \kappa + \epsilon]$. 

{\it Proof of Claim 5.}
Let
$$F(t,t',k) := \psi^b(t,k)(\psi^b)^{-1}(t',k) = \psi(t,k)\psi^{-1}(t',k).$$
Since $F(t,t',\cdot)$ is entire, $\lim_{k \to \kappa} F(t,t',k)$ exists and is bounded for fixed $t,t' \in [0, \infty)$. Moreover, for fixed $k \in (\kappa, \kappa + \epsilon]$, the identity (\ref{psibpsibidentity}) shows that $F(t,t',k)$ is bounded for $t, t' \in [0,\infty)$. On the other hand, the function $F(t,t',\kappa)$ is bounded for $t, t' \in [0,\infty)$ because $\psi^b(t,\kappa)$ is periodic or antiperiodic in $t$. 
This is enough to conclude that $\nu(t,k)$ is well-defined by (\ref{nuVolterra2}) for $(t,k) \in [0, \infty) \times [\kappa, \kappa + \epsilon]$. However, in order to prove the continuity we need estimates on $F$ as $t, t' \to \infty$ which are valid uniformly with respect to $k \in [\kappa, \kappa + \epsilon]$. 

We will prove the uniform bound
\begin{align}\label{Funiformbound}
|F(t,t',k)| \leq C(1 + |t-t'|), \qquad k \in [\kappa, \kappa + \epsilon], \quad t, t' \in [0, \infty).
\end{align}
For any integers $n,m \geq 0$, we have
\begin{align}\label{Fnm}
F(t + n\tau,t' + m\tau,k) = \psi(t,k) \mathcal{F}(k,m-n) \psi^{-1}(t',k),
\end{align}
where, for $j \in \Z$,
\begin{align*}
\mathcal{F}(k,j) & = S^b(k) e^{i\tilde{\Omega}(k)j\tau\sigma_3}S^b(k)^{-1}
	\\
& = 2i \sin(\tilde{\Omega}(k)j \tau) \begin{pmatrix} A^b(k)^2   & 
- A^b(k)B^b(k) 	\\
A^b(k) \overline{B^b(\bar{k})} & - A^b(k)^2 \end{pmatrix}
+ e^{-i\tilde{\Omega}(k)j\tau\sigma_3}.
\end{align*}
The functions $\sqrt{G}$, $z$, and $\tilde{\Omega}$ are analytic $U \to \C$, where $U$ denotes a small neighborhood of $\kappa$ in the Riemann surface (\ref{RiemannSurface}). Suppose $N \geq 1$ is the odd order of the zero of $G(k)$  at $\kappa$. Then, as an analytic function $U \to \C$ of $l = \sqrt{k - \kappa}$, $\sqrt{G}$ has a zero  of  order $N$ at $\kappa$. It follows that $z(k)$ equals $1$ or $-1$ to order $N$ at $\kappa$ as a function of $l$. Hence there exists an integer $p \in \Z$ such that $\tilde{\Omega}(k) - \frac{p\omega}{2}$ vanishes to order $N$ at $\kappa$ as an analytic function of $l$. 
It follows that the functions
\begin{align}\nonumber
& A^b(k)^2\Big(\tilde{\Omega}(k)- \frac{p\omega}{2}\Big) = -\frac{Z_{11} - Z_{22} - \sqrt{G}}{2\sqrt{G}}\Big(\tilde{\Omega}(k)- \frac{p\omega}{2}\Big), 
	\\ \nonumber
& A^b(k)B^b(k)\Big(\tilde{\Omega}(k)- \frac{p\omega}{2}\Big) = \frac{Z_{12}}{\sqrt{G}}\Big(\tilde{\Omega}(k)- \frac{p\omega}{2}\Big), 
	\\ \nonumber
& A^b(k)\overline{B^b(\bar{k})} \Big(\tilde{\Omega}(k)- \frac{p\omega}{2}\Big) = - \frac{Z_{21}}{\sqrt{G}}\Big(\tilde{\Omega}(k)- \frac{p\omega}{2}\Big),
\end{align}
are bounded for $k \in [\kappa, \kappa + \epsilon]$. 
Writing
\begin{align*}
\mathcal{F}(k,j) = &\; 2i \bigg(\frac{\sin(\tilde{\Omega}(k)j \tau)}{\tilde{\Omega}(k) j \tau - j p\pi}\bigg) j \tau \begin{pmatrix} A^b(k)^2   & 
- A^b(k)B^b(k) 	\\
A^b(k) \overline{B^b(\bar{k})} & - A^b(k)^2 \end{pmatrix}\Big(\tilde{\Omega}(k) - \frac{p\omega}{2}\Big) 
	\\
& + e^{-i\tilde{\Omega}(k)j\tau\sigma_3}.
\end{align*}
and using the bound
$$\bigg|\frac{\sin{x}}{x - j\pi}\bigg| \leq 1, \qquad x \in \R, \quad j \in \Z,$$
this implies
\begin{align}\label{mathcalFbound}
 |\mathcal{F}(k,j)| \leq C(1 + |j|), \qquad k \in [\kappa, \kappa + \epsilon], \quad j \in \Z.
\end{align}
Since $\psi(t,k)$ is bounded on the compact set $[0,\tau] \times [\kappa, \kappa + \epsilon]$, equations (\ref{Fnm}) and (\ref{mathcalFbound}) yield
$$|F(t + n\tau,t' + m\tau,k)| \leq C(1 + |m-n|), \qquad t,t' \in [0, \tau], \quad k \in [\kappa, \kappa + \epsilon], \quad m,n \geq 0.$$
The bound in (\ref{Funiformbound}) follows.

We write (\ref{nuVolterra2}) as
\begin{align}\nonumber
\tilde{\nu}(t,k) = &\; (1,0,0,1)^T - \int_t^\infty K(t,t',k) \tilde{\nu}(t', k) dt', \qquad t \in [0, \infty), \quad k \in [\kappa, \kappa + \epsilon],
\end{align}
where $\tilde{\nu} = (\nu_{11}, \nu_{12}, \nu_{21}, \nu_{22})^T$ and, in view of (\ref{u0uxasymptote}) and (\ref{Funiformbound}), the $4 \times 4$-matrix valued kernel $K(t,t',k)$ satisfies
$$|K(t,t',k)| \leq C (1 + t')^{-3/2}, \qquad k \in [\kappa, \kappa + \epsilon], \quad 0 \leq t \leq t' < \infty.$$
In particular,
$$\|K(t, \cdot, k)\|_{L^1([t, \infty))} \leq C (1+t)^{-1/2}, \qquad k \in [\kappa, \kappa + \epsilon], \quad 0 \leq t < \infty.$$
The same type of argument leading to (\ref{Psiminusf}) therefore implies that $\nu(t,k)$ satisfies 
$$|\nu(t,k) - I | \leq C (1+t)^{-1/2}, \qquad t \in  [0, \infty), \quad k \in [\kappa, \kappa + \epsilon],
$$
and is continuous on $[0,\infty) \times [\kappa, \kappa + \epsilon]$.
\proofendcontinue

\medskip\noindent
{\bf Claim 6.} For the defocusing NLS, $G(k)$ has no zeros of odd order in $\text{Int}(\bar{D}_1 \cup \bar{D}_4) \cap \R$.

{\it Proof of Claim 6.}   
The function $Q^b(k)$ is analytic from $U\setminus \{\kappa\}$ to $\C$ where $U$ is a small neighborhood of $\kappa$ in the Riemann surface (\ref{RiemannSurface}). Hence, we may write 
$$Q^b(k) = \sum_{j = -\infty}^\infty q_j l^{j}, \qquad k \in U \setminus \{\kappa\},$$
where $l = \sqrt{k - \kappa}$ and $\{q_j\}_{-\infty}^\infty$ are constants. By Claim 4, $|Q^b_\pm(k)| = 1$ for $k \in [\kappa - \epsilon, \kappa)$. Hence
$$Q^b(k) = \sum_{j = 0}^\infty q_j l^{j}, \qquad k \in U.$$
In particular, $Q^b(k)$ is bounded near $\kappa$ and $|Q^b(\kappa)| = 1$. 


By Claim 5, $\nu(0,k)$  is well-defined for $k \in [\kappa, \kappa + \epsilon]$. As in (\ref{fracBA}), we find
$$\frac{B(k)}{A(k)} = \frac{\nu_{11}(0,k) Q^b(k) + \nu_{12}(0,k)}{\nu_{21}(0,k)Q^b(k) + \nu_{22}(0,k)}, \qquad k \in [\kappa, \kappa + \epsilon].$$
In view of the global relation 
\begin{align*}
\frac{b(k)}{a(k)} = \frac{B_+(k)}{A_+(k)}, \qquad k \in \bar{D}_1 \cap \R,
\end{align*}
this gives
$$\bigg|\frac{b(\kappa)}{a(\kappa)}\bigg|
= \bigg|\frac{B(\kappa)}{A(\kappa)}\bigg|
= |Q^b(\kappa)| \bigg|\frac{\nu_{11}(0,\kappa)  + \nu_{12}(0,\kappa)\overline{Q^b(\kappa)}}{\nu_{21}(0,\kappa)Q^b(\kappa) + \nu_{22}(0,\kappa)}\bigg|.$$
But the symmetry $\nu(t, k) = \sigma_1 \overline{\nu(t, \bar{k})} \sigma_1$ which is valid for $k \in [\kappa, \kappa + \epsilon]$ then implies
$$\bigg|\frac{b(\kappa)}{a(\kappa)}\bigg|
= |Q^b(\kappa)| = 1.$$
Hence
$$\frac{1}{|a(\kappa)|^2} = 1 - \bigg|\frac{b(\kappa)}{a(\kappa)}\bigg|^2 = 0.$$
This contradicts the fact that $a(k)$ is bounded on $\R$.
\proofend

\begin{remark}\upshape
Under the assumptions of Theorem \ref{asymptoticallyadmissibleth}, we have
$$|Q^b_+(k)| \leq 1, \qquad k \in \bar{D}_1 \cap \R$$ 
for the defocusing NLS. Indeed, equation (\ref{mu1calE}) and the fact that $\mathcal{E}(n\tau, k) = \mathcal{E}(0,k) = S^b(k)$ for $n \in \Z$ imply that 
$$B(n\tau,k) \to B^b(k) \quad \text{and} \quad A(n\tau,k) \to A^b(k) \quad \text{as} \quad n \to \infty,$$
uniformly for $k$ in compact subsets of $\bar{D}_1 \cup \bar{D}_3$ that contain no branch points. Together with the global relation (\ref{gGR}) this implies the following pointwise convergence:
$$\frac{b(n\tau,k)}{a(n\tau,k)} = \frac{B_+(n\tau,k)}{A_+(n\tau,k)} \to \frac{B_+^b(k)}{A_+^b(k)} = Q^b_+(k), \qquad t \to \infty, \quad k \in (\bar{D}_1 \cap \R) \setminus \mathcal{D},$$
where $\mathcal{D}$ denotes the set of branch points together with the set of zeros of $A_+^b(k)$.
For the defocusing NLS, $\frac{b(t,k)}{a(t,k)}$ is a smooth function of $k \in \R$ which satisfies
$$\bigg|\frac{b(t,k)}{a(t,k)}\bigg|^2 = 1 - \frac{1}{|a(t,k)|^2} < 1, \qquad k \in \R,$$
for each $t \geq 0$. Since $\mathcal{D}$ is a discrete set, it follows that $|Q^b_+(k)| \leq 1$ for $k \in \bar{D}_1 \cap \R$.
\end{remark}

We next illustrate the general formalism with some examples for which the background $t$-part can be solved explicitly.

\section{Example: Stationary one-solitons}\label{stationarysec}\nequation
The focusing NLS admits the following family of stationary soliton solutions:
\begin{align}\label{stationarysoliton}
u(x,t) = \frac{\sqrt{\omega}}{\cosh(x\sqrt{\omega} - \gamma)}  e^{it\omega}, \qquad \gamma \in \R, \quad \omega > 0.
\end{align}
These are smooth solutions which decay as $x \to \infty$.
By evaluating $u(x,t)$ and $u_x(x,t)$ at $x = 0$, we find the following Dirichlet and Neumann boundary values: 
$$g_0(t) = g_0^b(t) = \alpha e^{i\omega t}, \qquad g_1(t) = g_1^b(t) = ce^{i\omega t},$$
where the constants $\alpha$ and $c$ are given by
\begin{align}\label{stationaryac}
\alpha = \frac{\sqrt{\omega}}{\cosh(\gamma)} > 0, \qquad c = \omega \frac{\sinh(\gamma)}{\cosh^2(\gamma)} \in \R.
\end{align}
The functions $g_0(t)$ and $g_1(t)$  are periodic functions with period $\tau = \frac{2\pi}{\omega}$.
We write
\begin{align}\label{stationaryc}
c = \sigma \alpha\sqrt{\omega - \alpha^2} \quad \text{with} \quad \sigma = \sgn(\gamma).
\end{align}
Direct integration of the $x$-part of (\ref{lax}) yields
$$\mu_3(x,t,k) = \begin{pmatrix}  \frac{2 k+i \sqrt{\omega } \tanh (\gamma -x \sqrt{\omega })}{2 k-i
   \sqrt{\omega }} & -\frac{e^{i t \omega } \sqrt{\omega }
   \sech(\gamma -x \sqrt{\omega })}{\sqrt{\omega }-2 i k} \\
 \frac{e^{-i t \omega } \sqrt{\omega } \sech(\gamma -x \sqrt{\omega
   })}{2 i k+\sqrt{\omega }} & \frac{2 k-i \sqrt{\omega } \tanh
   (\gamma -x \sqrt{\omega })}{2 k+i \sqrt{\omega }}\end{pmatrix}.$$
Evaluation of this expression at $x = t = 0$ gives
\begin{align}\label{stationaryab}
b(k) = \frac{\sqrt{\omega} \sech(\gamma)}{2ik - \sqrt{\omega}}, \qquad
 a(k) = \frac{2k - i\sqrt{\omega}\tanh(\gamma)}{2k + i\sqrt{\omega}}.
\end{align} 
The relation $\mu_2 = \mu_3 e^{-i(kx + 2k^2 t)\hat{\sigma}_3}s^{-1}$ yields
\begin{align*}
(\mu_2(x,t,k))_{12} = &\; \frac{\sqrt{\omega }}{4 k^2+\omega }  \Big[(\sqrt{\omega } \tanh
   (\gamma )-2 i k) e^{i t \omega} \sech(\gamma -x \sqrt{\omega })
   	\\
&   + e^{-2i(k x + 2k^2 t)} \sech(\gamma ) (2 ik -\sqrt{\omega } \tanh (\gamma -x \sqrt{\omega }))\Big],
   	\\
(\mu_2(x,t,k))_{22} = &\; \frac{1}{4 k^2+\omega }\Big[\omega  \sech(\gamma ) e^{-i t\omega} e^{-2i(kx + 2k^2 t)}
   \sech(\gamma -x \sqrt{\omega })
   	\\
& + (2 k+i \sqrt{\omega } \tanh (\gamma )) (2 k-i \sqrt{\omega } \tanh (\gamma -x
   \sqrt{\omega }))\Big].
\end{align*}
Using the relation $Z(k) = \mu_2(0,\tau,k) e^{-2ik^2 \tau\sigma_3}$ together with (\ref{lambdaGdef}) and (\ref{branchchoices}), we find
\begin{align}\nonumber
& \sqrt{G(k)} = 2i\sin(2k^2 \tau), \qquad z(k) = e^{-2ik^2\tau}, \qquad \tilde{\Omega}(k) = 2k^2,
	\\ \label{stationaryE}
& (\mathcal{E}(t,k))_{12} = \frac{\sqrt{\omega } e^{i t \omega } \sqrt{1-\frac{\omega  \text{sech}^2(\gamma )}{4 k^2+\omega }}}{\sqrt{\omega } \sinh (\gamma )+2 i k \cosh (\gamma )}, \qquad
(\mathcal{E}(t,k))_{22} = \sqrt{1-\frac{\omega  \text{sech}^2(\gamma )}{4 k^2+\omega }}.
\end{align}	
\begin{figure}
\begin{center}
\begin{overpic}[width=.55\textwidth]{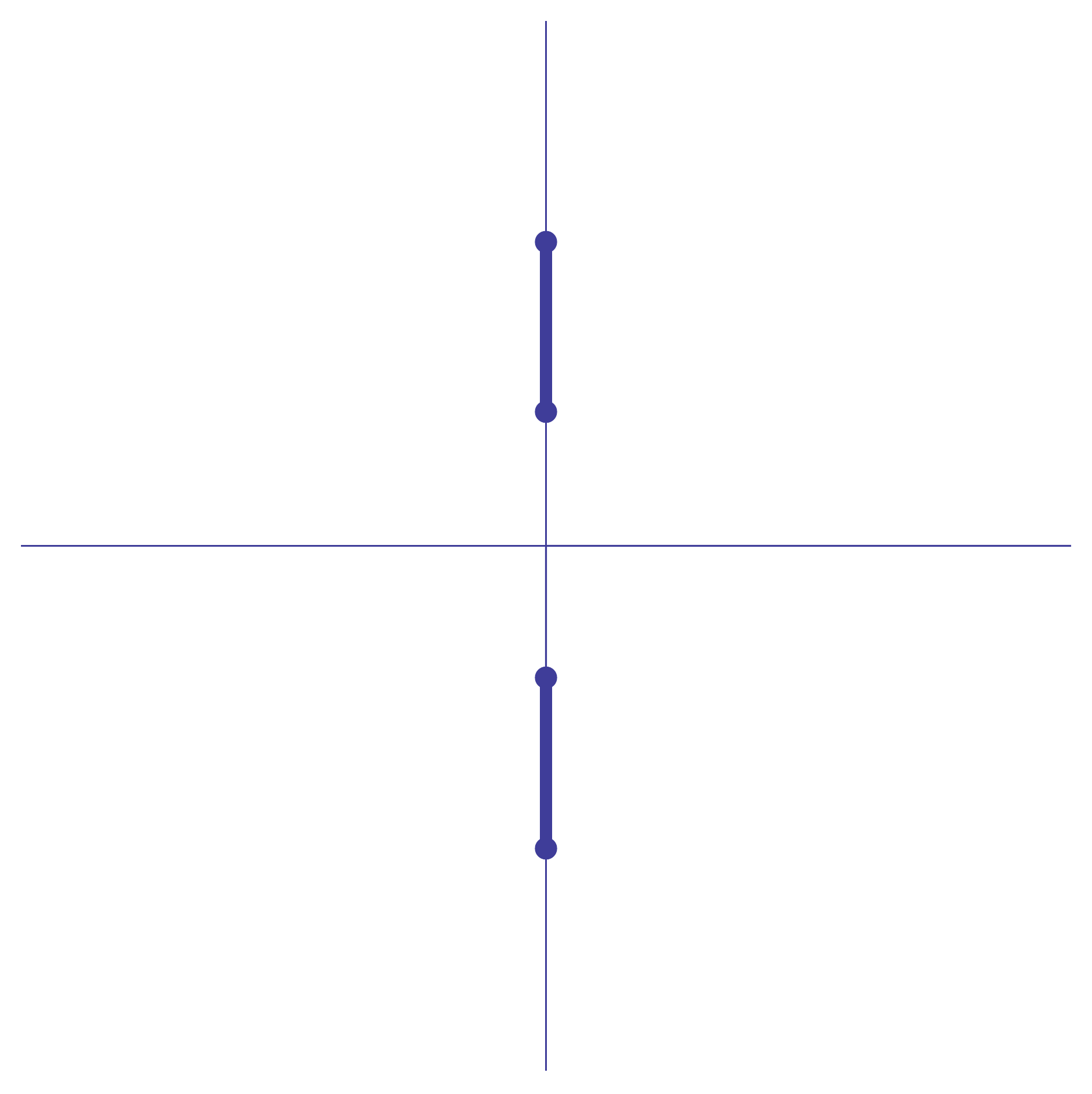}
     \put(80,80){$D_1$}
     \put(16,80){$D_2$}
     \put(16,20){$D_3$}
     \put(80,20){$D_4$}
      \put(52,77){$\frac{i\sqrt{\omega}}{2}$}
      \put(52,61){$\frac{i\sqrt{\omega}\tanh|\gamma|}{2}$}
      \put(52,37){$-\frac{i\sqrt{\omega}\tanh|\gamma|}{2}$}
      \put(52,21.5){$-\frac{i\sqrt{\omega}}{2}$}
\end{overpic}
     \begin{figuretext}\label{stationarycuts.pdf}
         The domains $\{D_j\}_1^4$ and the branch cuts $\mathcal{C}$ in the case of the stationary one-soliton (\ref{stationarysoliton}).
      \end{figuretext}
     \end{center}
\end{figure}
Hence the $D_j$'s coincide with the four quadrants of the complex $k$-plane and all functions are well-defined in $\C \setminus \mathcal{C}$ where $\mathcal{C} \subset i\R$ consists of two cuts along the imaginary axis (see Figure \ref{stationarycuts.pdf})
$$\mathcal{C} = \bigg[-\frac{i\sqrt{\omega}}{2}, - \frac{i\sqrt{\omega}\tanh |\gamma|}{2}\bigg] \cup \bigg[\frac{i\sqrt{\omega}\tanh|\gamma|}{2}, \frac{i\sqrt{\omega}}{2}\bigg].$$
We find $S(k)$ from equation (\ref{sS}), which then gives
\begin{align}\label{stationaryAB}
B(k) = \frac{\sqrt{\omega } \sqrt{1-\frac{\omega  \sech^2(\gamma )}{4 k^2+\omega}}}{\sqrt{\omega } \sinh (\gamma )+2 i k \cosh (\gamma )}, \qquad
A(k) = \sqrt{1-\frac{\omega  \sech^2(\gamma )}{4 k^2+\omega }}.
\end{align}
Finally, the identity $\mu_1 = \mu_2 e^{-ikx\hat{\sigma}_3}\big[e^{-2ik^2t\sigma_3}S(k)e^{i\tilde{\Omega}(k) t\sigma_3}\big]$ gives
\begin{align} \nonumber
(\mu_1(x,t,k))_{12} = &\; A(k) \frac{\sqrt{\omega } e^{i t \omega } \sech(\gamma -x \sqrt{\omega
   })}{\sqrt{\omega } \tanh (\gamma )+2 i k},
   	\\ \label{stationarymu1}
(\mu_1(x,t,k))_{22} = &\; A(k) \frac{\sqrt{\omega } \tanh (\gamma -x \sqrt{\omega })+2 i
   k}{\sqrt{\omega } \tanh (\gamma )+2 i k}.
\end{align}

It is easy to verify that the spectral functions given in (\ref{stationaryab}) and (\ref{stationaryAB}) satisfy the global relation:
$$b(k)A(k) - a(k) B(k) = 0, \qquad k \in D_1.$$
In fact, denoting the branch points by $\pm K_1$ and $\pm K_2$ where 
$$K_1 = \frac{i \sqrt{\omega}}{2}, \qquad K_2 = \frac{i\bar{c}}{2a} = \frac{i \sqrt{\omega}}{2} \tanh (\gamma ),$$ 
we have
$$ b(k) = \frac{\alpha}{2i(k + K_1)}, \qquad a(k) = \frac{k - K_2}{k + K_1},$$
and
$$B(k) = \frac{\alpha}{2i}\sqrt{\frac{k + K_2}{(k + K_1)(k-K_1)(k-K_2)}}, \qquad 
A(k) = \sqrt{\frac{(k + K_2)(k-K_2)}{(k + K_1)(k - K_1)}}. $$

Clearly, $G(k)$ has no zeros of odd order. In fact, $G(k) = -4\sin^2(2k^2\tau)$ has double zeros at each point in the set
$$\Big\{\pm \frac{\sqrt{n\omega}}{2}, \pm \frac{i\sqrt{n\omega}}{2} \, \Big| \, n = 1,2,\dots\Big\}$$
and a zero of fourth order at the origin. 

\begin{remark}\upshape
Equation (\ref{stationaryc}) shows that $\{\alpha e^{i\omega t}, ce^{i\omega t}\}$ is an asymptotically admissible pair for the focusing NLS provided that the triple $(\alpha,\omega,c)$ satisfies the condition (\ref{1.3a}).
\end{remark}

\section{Example: Single exponentials}\upshape
Consider the case when both $g_0^b$ and $g_1^b$ are single exponentials of the same frequency:
\begin{align}\label{singleexp}
g_0^b(t) = \alpha e^{i\omega t}, \qquad g_1^b(t) = c e^{i\omega t}, \qquad \alpha > 0, \quad c \in \C, \quad \omega \in \R.
\end{align}
The special case when $\alpha$ and $c$ are given by (\ref{stationaryac}) was studied in section \ref{stationarysec}.
However, the background $t$-part can be solved explicitly for any values of $\alpha$ and $c$.
Indeed, direct integration of the background $t$-part (\ref{backgroundtpart}) yields
\begin{align*}
  & \psi(t,k) = e^{\frac{i\omega}{2} t \sigma_3} \begin{pmatrix} \cos(t\Omega) + \frac{4k^2 + 2\lambda \alpha^2 + \omega}{2i \Omega}\sin(t\Omega) & (2\alpha k + ic)\frac{\sin(t\Omega)}{\Omega} \\
 \lambda (2\alpha k  - i\bar{c})\frac{\sin(t\Omega)}{\Omega} & \cos(t\Omega) - \frac{4k^2 + 2\lambda \alpha^2 + \omega}{2i\Omega}\sin(t\Omega)
 \end{pmatrix},
\end{align*}
where $\Omega(k)$ is defined by
\begin{align} \label{Omega2def}
& \Omega(k) = \sqrt{4k^4 + 2\omega k^2 + 4\lambda \alpha \im(c) k + \bigg(\frac{\omega}{2} + \lambda \alpha^2\bigg)^2 - \lambda |c|^2}.
\end{align}
We fix the branch of $\Omega(k)$ by requiring that
$$\Omega(k) = 2k^2 + \frac{\omega}{2} + O(k^{-1}), \qquad k \to \infty.$$
Letting $\tau = 2\pi/\omega$ and choosing the branches according to (\ref{branchchoices}), we find
\begin{align*}
& \tr{Z(k)} = -2 \cos(\Omega(k) \tau), \qquad Z_{12}(k) = -(2\alpha k + ic) \frac{\sin(\Omega(k) \tau)}{\Omega(k)},
	\\
& Z_{11}(k) - Z_{22}(k) = i(4k^2 + 2\lambda \alpha^2 + \omega) \frac{\sin(\Omega(k) \tau)}{\Omega(k)},
	\\
& z(k) = -e^{-i \Omega(k) \tau}, \qquad \tilde{\Omega}(k) = \Omega(k) - \frac{\omega}{2}, \qquad 
G(k) = -4\sin^2(\Omega(k)\tau), 
	\\
& \sqrt{G(k)} = 2i\sin(\Omega(k) \tau).
\end{align*}
Using the identities
\begin{align*}
&\frac{Z_{11} - Z_{22} - \sqrt{G}}{\sqrt{G}} = \frac{H - 2\Omega}{\Omega}, \qquad
 \frac{Z_{11} - Z_{22} + \sqrt{G}}{\sqrt{G}} = \frac{H}{\Omega},
 	\\
&  (H-2\Omega)H =\lambda (2\alpha k - i\bar{c})(2\alpha k + ic),
\end{align*}
where $H(k)$ is defined by
\begin{align}\label{Hdef}
H(k) = \Omega(k) - 2k^2 - \lambda \alpha^2 - \frac{\omega}{2},
\end{align}
we infer the following formulas:
\begin{align*}
& S^b(k) = \sqrt{\frac{2\Omega - H}{2\Omega}} 
\begin{pmatrix} 1 & \frac{c - 2i\alpha k}{2\Omega - H} \\
\lambda \frac{\bar{c} + 2i\alpha k}{2\Omega - H} & 1 \end{pmatrix} 
= \sqrt{\frac{2\Omega - H}{2\Omega}} \begin{pmatrix} 1 & \frac{\lambda iH}{2\alpha k - i\bar{c}} \\
-\frac{iH}{2\alpha k + ic} & 1 \end{pmatrix},
	\\
& e^{t\mathcal{B}(k)} = \begin{pmatrix}
 \cos (\tilde{\Omega} t) - i \frac{\Omega -H}{\Omega} \sin (\tilde{\Omega} t) &
 \frac{2\alpha k + ic}{\Omega}  \sin (\tilde{\Omega} t) \\
 \frac{2\alpha k - i\bar{c}}{\Omega} \lambda  \sin (\tilde{\Omega} t) &  \cos (\tilde{\Omega} t) + i\frac{\Omega - H}{\Omega}  \sin (\tilde{\Omega} t)
   \end{pmatrix},
   	\\
& P(t,k) = \frac{1}{2 \Omega} \begin{pmatrix} e^{it\omega} H + 2\Omega - H 
& (e^{i t \omega} - 1) (c-2 i \alpha k)  \\
\lambda (e^{-it\omega} - 1)(\bar{c} + 2i\alpha k) & He^{-it\omega} + 2\Omega - H \end{pmatrix},
	\\
& \mathcal{E}(t,k) = e^{\frac{i\omega}{2}t\hat{\sigma}_3}S^b(k). 		
\end{align*}
The background eigenfunction $\psi^b$ defined in (\ref{psiB}) is given by
\begin{align}\label{psiBexample}
\psi^b(t,k) = e^{\frac{i\omega}{2}t \sigma_3} S^b(k) e^{-i \Omega(k)  t \sigma_3}.
\end{align}
The set of branch points is given by
$$\mathcal{P} = \bigg\{-\frac{ic}{2\alpha}, \frac{i\bar{c}}{2\alpha}, \text{zeros of }  \Omega^2\bigg\}.$$

\begin{remark}\upshape
In the case of $\lambda = -1$, the eigenfunction $\psi^b(t,k)$ in (\ref{psiBexample}) coincides with the background eigenfunction adopted in \cite{BKS2009}.
\end{remark}

Writing
$$G(k) = -4 \bigg(\frac{\sin(\Omega(k) \tau)}{\Omega(k)}\bigg)^2 \Omega^2(k)$$
and noting that $\frac{\sin(\Omega(k) \tau)}{\Omega(k)}$ is an entire function, we see that the zeros of odd order of $G(k)$ are exactly the zeros of odd order of $\Omega^2(k)$.
In the special case of single exponential profiles, Theorem \ref{asymptoticallyadmissibleth} therefore reduces to the following result.

\begin{theorem}\label{asymptoticallyadmissibleth2}
 Let $\{g_0^b(t), g_1^b(t)\}$ be given by (\ref{singleexp}). If the pair $\{g_0^b(t), g_1^b(t)\}$ is asymptotically admissible for NLS, then the function $\Omega^2(k)$, defined in (\ref{Omega2def}), has no zeros of odd order in the interior of $\bar{D}_1 \cup \bar{D}_4$.
\end{theorem}

For the focusing NLS, the imposition of the requirement of Theorem \ref{asymptoticallyadmissibleth2} leads to the two families of  parameter triples in (\ref{1.3}), see \cite{BKS2009}. 

For the defocusing NLS, imposing the analogous requirement leads to {\it five} different families of triples; see \cite{Ldefocusing} for the relatively lengthy proof of the following result. 

\begin{theorem}\label{mainth}
If the pair $\{\alpha e^{i\omega t}, ce^{i\omega t}\}$ is asymptotically admissible for the defocusing NLS equation, where $\alpha>0$, $\omega \in \R$, and $c \in \C$, then the triple $(\alpha, \omega, c)$ belongs to one of the following disjoint subsets:
\begin{align*}
& \bigg\{\bigg(\alpha, \omega, c = \pm \sqrt{\frac{(\omega + 3\alpha^2)^3}{27\alpha^2}} + \frac{i|\omega|^{3/2}}{3\sqrt{3} \alpha}\bigg) \; \bigg| \;  \alpha > 0, \; -3\alpha^2 \leq \omega < 0\bigg\},
	\\ \nonumber
& \bigg\{\bigg(\alpha  = -\frac{4K^3 + \omega K}{c_2}, \omega, c = \pm \sqrt{ \bigg(\alpha^2+ \frac{\omega}{2}\bigg)^2 -c_2^2 - 2K^2(6K^2 + \omega)} + ic_2\bigg)  
	\\ 
&\hspace{3.5cm} \bigg| \;  -12 K^2 < \omega < -4K^2, \; 0 < c_2 \leq -\frac{4K^2 + \omega}{2}, \; K > 0\bigg\}, 
	\\ 
& \big\{(\alpha,\omega, c = i\alpha\sqrt{-2\alpha^2 - \omega}) \; \big| \; \alpha > 0, \; \omega < -3\alpha^2\big\},
	\\ 
& \big\{(\alpha,\omega,c = \pm \alpha\sqrt{\omega + \alpha^2}) \; \big| \; \omega + \alpha^2 \geq 0, \; \alpha > 0\big\},
	\\ \nonumber
& \bigg\{\bigg(\alpha = -\frac{4K^3 + \omega K}{c_2}, \omega, c = \pm \sqrt{\bigg(\alpha^2+ \frac{\omega}{2}\bigg)^2 -c_2^2 - 2K^2(6K^2 + \omega)} + ic_2\bigg)
	\\ 
& \hspace{3.5cm} \bigg| \; -4K^2 < \omega \leq -3K^2, \; -\frac{4K^2 + \omega}{2} \leq c_2 < 0, \; K > 0\bigg\}.
\end{align*}
\end{theorem}

In addition to the analogs of the branches (\ref{1.3}) present in the focusing case, in the defocusing case there are three branches for which both the real and the imaginary parts of $c$ are nonzero. Thus, even for the relatively simple example of a single exponential, the analysis of the defocusing NLS is surprisingly rich. Since each of the additional branches depends on two or three parameters, this provides a large number of potentially asymptotically admissible pairs for the defocusing NLS.

\section{Asymptotically consistent pairs}\nequation\label{consistentsec}
Theorem \ref{asymptoticallyadmissibleth} shows that the condition that $G(k)$ has no zeros of odd order in the interior of $\bar{D}_1 \cup \bar{D}_4$ is a necessary requirement for a pair to be asymptotically admissible. However, this condition may not be sufficient. We therefore make the following definition.

\begin{definition}
We call $\{g_0^b(t), g_1^b(t)\}$ an {\it asymptotically consistent pair} if $G(k)$ has no zeros of odd order in the interior of $\bar{D}_1 \cup \bar{D}_4$ for any choice\footnote{The choice is subject to the underlying assumptions that  $\mathcal{C}$  is invariant under $k \to \bar{k}$  and $D_1\setminus \mathcal{C}$ is connected.}
of branch cuts. 
\end{definition}

We propose two different approaches for identifying among the asymptotically consistent pairs those that are actually asymptotically admissible. The first approach involves summing up a perturbative series and is presented in detail in \cite{tperiodicII}. The second approach is a generalization of the approach used for the single exponential profile (\ref{1.1}) for the focusing NLS in \cite{BK2007, BKS2009}. It involves the study of the RH problem with jump contour
$$\{k \in \C \, | \, \im k = 0 \;\, \text{or}\, \im \tilde{\Omega}(k) = 0\} \cup \mathcal{C},$$
satisfied by the following sectionally analytic function:
\begin{align*}
  M(x,t,k) = \begin{cases}
  \begin{pmatrix} \frac{[\phi_2(x,t,k)]_1}{a(k)} & [\phi_3(x,t,k)]_2 \end{pmatrix} e^{i(kx + \tilde{\Omega}(k) t)\sigma_3}, & k \in D_1  \setminus \mathcal{C},
   \vspace{.1cm}	 \\ 
  \begin{pmatrix} \frac{[\phi_1(x,t,k)]_1}{d(k)} & [\phi_3(x,t,k)]_2 \end{pmatrix} e^{i(kx + \tilde{\Omega}(k) t)\sigma_3}, & k \in D_2 \setminus \mathcal{C},
   \vspace{.1cm}	\\ 
  \begin{pmatrix} [\phi_3(x,t,k)]_1 & \frac{[\phi_1(x,t,k)]_2}{\overline{d(\bar{k})}} \end{pmatrix} e^{i(kx + \tilde{\Omega}(k) t)\sigma_3}, & k \in D_3 \setminus \mathcal{C},
   \vspace{.1cm}	\\ 
  \begin{pmatrix} [\phi_3(x,t,k)]_1 & \frac{[\phi_2(x,t,k)]_2}{\overline{a(\bar{k})}} \end{pmatrix} e^{i(kx + \tilde{\Omega}(k) t)\sigma_3}, & k \in D_4  \setminus \mathcal{C},
\end{cases}	
\end{align*}
where $d(k) = a(k)\overline{A(\bar{k})} -  \lambda b(k) \overline{B(\bar{k})}$.
Asymptotic admissibility can be established by (a) proving the unique solvability of this RH problem, (b) defining the function $u(x,t)$ by
$$u(x,t) = 2i \lim_{k \to \infty} (kM(x,t,k))_{12},$$
and (c) showing that $u(x,t)$ solves the NLS equation (\ref{nls}) and satisfies (\ref{u0uxasymptote}).
In the case of the single exponential profile (\ref{1.1}) for the focusing NLS, asymptotic admissibility was established in this way in \cite{BK2007, BKS2009}.

\section{Concluding remarks}\nequation\label{concludesec}
In the series of papers \cite{BIK2007, BIK2009, BKS2009, BKSZ2010}, Boutet de Monvel and coauthors analyzed the focusing NLS equation with boundary conditions given asymptotically by a pair of single exponentials:
$$u(0,t) \sim \alpha e^{i\omega t}, \qquad u_x(0,t) \sim c e^{i\omega t}, \qquad t \to \infty,$$
where $a > 0$, $c \in \C$, and $\omega \in \R$.
By imposing a requirement derived from the global relation (essentially the requirement of Theorem \ref{asymptoticallyadmissibleth2} above), they were able to express $c$ in terms of $\alpha$ and $\omega$. Their result determines which of the single exponential pairs are {\it asymptotically admissible} and can be seen as a (partial) characterization of the Dirichlet to Neumann map for single exponential profiles in the limit of large $t$. 
In \cite{Ldefocusing}, single exponential profiles for the {\it defocusing} NLS were analyzed in an analogous way, again using the requirement of Theorem \ref{asymptoticallyadmissibleth2}. 

In this paper, in an effort to generalize the results of \cite{BIK2007, BIK2009, BKS2009, BKSZ2010, Ldefocusing}, we have proved a result (Theorem \ref{asymptoticallyadmissibleth}) which is analogous to Theorem \ref{asymptoticallyadmissibleth2} but which applies to any periodic pair of Dirichlet and Neumann profiles. In the case of single exponentials, Theorem \ref{asymptoticallyadmissibleth} reduces to Theorem \ref{asymptoticallyadmissibleth2}.

In general, the nonexplicit nature of the associated eigenfunction $\psi^b$ makes it challenging to apply Theorem \ref{asymptoticallyadmissibleth}. However, there exists a large collection of pairs for which $\psi^b$ is known explicitly. This collection includes the class of stationary solitons, as well as the class of finite-gap solutions generated by the Baker-Akhiezer formalism.
We propose that a large class of asymptotically consistent pairs can be generated by taking periodic finite-gap solutions $\psi^b$ of the background $t$-part and then applying the requirement of Theorem \ref{asymptoticallyadmissibleth}.

We emphasize that Theorems \ref{asymptoticallyadmissibleth} and \ref{asymptoticallyadmissibleth2} provide {\it necessary} conditions for a pair to be asymptotically admissible. In the special case of single exponential profiles for the focusing NLS, the results of \cite{BIK2007, BIK2009, BKS2009, BKSZ2010} show that these conditions are also sufficient, but this may not be the case in general. We refer to pairs fulfilling the requirement of Theorem \ref{asymptoticallyadmissibleth} as {\it asymptotically consistent}. We have proposed two different approaches for identifying among the asymptotically consistent pairs those that are actually asymptotically admissible. The first approach involves summing up a perturbative series and is presented in detail in \cite{tperiodicII}. The second approach involves the study of a RH problem and generalizes the approach used in \cite{BK2007, BKS2009} for the analysis of a single exponential for the focusing NLS.

\bigskip
\noindent
{\bf Ethics statement.} This work did not involve any collection of human data.

\medskip
\noindent
{\bf Data accessibility.} This work does not rely on any experimental data.

\medskip
\noindent
{\bf Competing interests.} We have no competing interests.

\medskip
\noindent
{\bf Authors' contributions.} JL and ASF conceived the ideas and proved the mathematical results presented in this paper. Both authors gave final approval for publication.

\medskip
\noindent
{\bf Acknowledgements.} The authors are grateful to the two referees for many helpful suggestions.

\medskip
\noindent
{\bf Funding statement.} This work was supported by the EPSRC, grant EP/H04261X/1. It was also co-financed by the European Union
(European Social Fund ESF) and Greek national funds through the
Operational Program Education and Lifelong Learning of the National
Strategic Reference Framework (NSRF) - Research Funding Program: THALIS
(MIS 379416). Investing in knowledge society through the European Social
Fund.

\bibliographystyle{plain}
\bibliography{is}

\end{document}